\documentclass[a4paper,12pt]{article}
\font\header=cmssdc10 at 20pt

\usepackage[ansinew]{inputenc}   
\usepackage{amsfonts}
\usepackage{indentfirst}
  \usepackage[normalem]{ulem}
  \usepackage{graphicx}
  \usepackage{amssymb}
 \usepackage{draftcopy}
  \usepackage{amsmath}
  \usepackage{rotating}
  \usepackage{color}
  \usepackage{setspace}
  \usepackage{fancybox}
 \usepackage{fancyvrb}
	\usepackage{fancyhdr}
	\usepackage[top=0.5cm,left=1cm,right=1cm,bottom=0.5cm]{geometry}
  \usepackage{wrapfig}
    \usepackage{float}
  \usepackage{type1cm}
\usepackage{eso-pic}
\usepackage{color}
 \usepackage{amsthm,amsmath,amsfonts}
\usepackage{fullpage,enumerate}
\usepackage{tikz}
\usepackage[colorlinks,citecolor=blue,urlcolor=blue]{hyperref}
\usepackage{hypernat}
\usepackage[notref,notcite]{showkeys}
\usepackage{graphicx}

\def\1{{\bf 1}}

\begin{document}

{\header Randomness versus selection in genome evolution}

\vskip1cm

Rinaldo B. Schinazi

 Department of Mathematics, University of Colorado,
 
 Colorado Springs, CO 80933-7150, USA; e-mail: 
Rinaldo.Schinazi@uccs.edu

\vskip1cm

{\bf Abstract.} We propose a Markov chain approach for the evolution of a genealogical line of genomes. Our idealized genome has $N$ sites and each site can be in state $0$ or $1$.
At each time step we pick a site at random. If the site is in state $0$ we flip it to state 1 with probability $p$ or we keep it in state $0$ with probability $1-p$. If the site is in state $1$ we flip it to state 0 with probability $1-p$ or we keep it in state $1$ with probability $p$.  Even when state 1 has a selective advantage (i.e. $p>1/2$) the Markov chain is quite unlikely to approach the most fit allele (i.e. all 1's). In fact, randomness (i.e. which site is picked for a possible mutation) and selection (i.e. the value of $p$) balance each other out so that the number of $1$'s in the genome converges to a Gaussian distribution centered around $Np$.

\vskip1cm


\vskip1cm

{\header  A generalized Ehrenfest model}

\vskip1cm

Consider a genome with $N$ sites. Each site can be in state $1$ or $0$. 
We pick a site uniformly. That is, if $k$ sites are in state $1$  the probability of picking a site in state $1$ is $k/N$ and the probability of picking a site in state $0$ is $(N-k)/N$.

\begin{itemize}
\item If we pick a site in state $0$ we flip it to state 1 with probability $p$ or we keep it in state $0$ with probability $1-p$. 
\item If we pick a site in state $1$ we flip it to state 0 with probability $1-p$ or we keep it in state $1$ with probability $p$. 
\end{itemize}

Let $Y_n$ be the number of 1's at time $n\geq 0$. Then, $(Y_n)$ is a finite discrete time Markov chain. At time $n\geq 0$, $Y_n=k$ for some positive integer $k$ between $0$ and $N$. Then at time $n+1$ there are three possibilities if $1\leq k\leq N-1$.
\begin{itemize}
\item With probability $p(k,k+1)=p\frac{N-k}{N}$ there is a mutation from $0$ to $1$ and therefore $Y_{n+1}=k+1$.
\item With probability $p(k,k-1)=(1-p)\frac{k}{N}$ there is a mutation from $1$ to $0$ and therefore $Y_{n+1}=k-1$.
\item With probability $1-p(k,k+1)-p(k,k-1)$ nothing happens and  $Y_{n+1}=k$.
\end{itemize}
Note that we have reflecting barriers at $0$ and $N$. That is, if $Y_n=0$ then $Y_{n+1}=1$ with probability $p$ or $Y_{n+1}=0$ with probability $1-p$. If $Y_n=N$ then $Y_{n+1}=N$ with probability $p$ or $Y_{n+1}=N-1$ with probability $1-p$.

The Markov chain $(Y_n)$ is finite and irreducible and therefore has a unique stationary distribution $\nu_N$. That is, if we pick $Y_0$ according to the probability distribution $\nu_N$ then at any given time $n$,  $\nu_N$ is the probability distribution of $Y_n$. Moreover, from any initial state the distribution of the chain converges exponentially fast to $\nu_N$,  see Bhattacharya and Waymire (1990) for this and other results for finite Markov chains.

Actually, $\nu_N$ is reversible, a  property that implies stationarity. That is, for every integer $k$ between $0$ and $N$,
$$\nu_N(k)p(k,k+1)=\nu_N(k+1)p(k+1,k).$$
Using these equations and the fact that $\sum_{k=0}^N\nu_N(k)=1$ we get for $0\leq k\leq N$,
$$\nu_N(k)=\begin{pmatrix}N\\k\end{pmatrix} p^k(1-p)^{N-k}.$$
That is, $\nu_N$ is a binomial distribution with parameters $N$ and $p$. By the Central Limit Theorem as $N\to+\infty$, $\nu_N$ will approach a Gaussian distribution centered at state $pN$. 

Even when state 1 has a selective advantage (i.e. $p>1/2$) the Markov chain is quite unlikely to approach the most fit allele (i.e. all 1's). In fact, randomness (i.e. which site is picked for a possible mutation) and selection (i.e. the value of $p$) balance each other out so that the number of $1$'s in the genome converges to a Gaussian distribution centered around $Np$.
Intuitively, we see the following picture.  As the number of 1's increases, a site in state 1 is more likely than a site in state 0 to get picked at random. A $1-p$ fraction of these 1's will flip to 0. This prevents the number of 1's to continue increasing to $N$. Instead a (random) equilibrium for the number of 1's establishes itself around $Np$. This is confirmed by solving the equation $p(k,k+1)=p(k,k-1)$ which yields $k=Np$.

\vskip1cm

{\header  Time to to reach the most fit allele}

\vskip1cm

From the stationary distribution $\nu_N$ we see that for all values $p<1$ the chain will very likely stay away from the most fit allele. Worse than that, even if we start the chain $(Y_n)$ in state $N$ then 
the expected time for the chain to return to state $N$ is known to be
$$\frac{1}{\nu_N(N)}=(\frac{1}{p})^N.$$
That is, as $N\to+\infty$, the expected return time grows exponentially fast. It is interesting to compare this time to the expected time it takes to return to state $pN$. Assuming $pN$ is an integer, by Stirling's formula as $N\to +\infty$ we have
$$\frac{1}{\nu_N(Np)}\sim \left(2\pi p(1-p)\right)^{1/2} N^{1/2}.$$
Hence, the expected time to return to $pN$ is only of order $N^{1/2}$.

Observe that the selection parameter $p$ plays only a rather marginal role in the model (except when $p=1$, see below). Having $p<1/2$ or $p>1/2$ does not change the qualitative behavior of the model. The parameter $p$ plays only a role in that the most likely alleles are around $pN$.

\bigskip

$\bullet$ {\bf The case p=1.} In this case a $0$ can flip to $1$ but a $1$ cannot flip. That is, the genome can only become more and more fit. This corresponds to a 'directional' evolution (which we know does not occur!). Only in this case is the genome converging to the all 1's allele as we show below.

The transition probabilities for the chain $(Y_n)$ become for $p=1$ and  $0\leq k\leq N$,
\begin{itemize}
\item  $p(k,k+1)=\frac{N-k}{N}$.
\item  $p(k,k)=\frac{k}{N}$ 
\end{itemize}
In particular, the all $1$'s allele is absorbing. 
The finite Markov chain $(Y_n)$ will eventually get absorbed by the state $N$. The time it takes to get absorbed is the well-known collector's problem, see for instance Port (1994). Starting from all $0$'s the expected time to reach all $1$'s is of order $N\ln N$ as $N$ goes to infinity. Hence, in this particular case the most fit allele is attained and this happens in a relatively short time.

\vskip1cm

{\header Literature}

\vskip1cm

There is a vast literature using probability models for the evolution of a genome.  Most models require
every step to increase the fitness of the genome, see Gillepsie (1983), Kaufman and Levin (1987), Orr (2003), Hegarty and Martinsson (2014) and Schinazi (2019). These models can be thought of as population models that keep track of the most fit genome in the population at every time step. That is, they do not follow a genealogical line, instead at every time step the model tracks the genome with the maximum fitness in the whole population. This is why the fitness is constrained to increase or stay put at every step. 

Our point of view is different. We model the evolution of one genealogical line of genomes. At every time step 
a new individual is born. Its genome is exactly the same as its mother's except possibly at a single site. We show that this sequence of genomes converges to a stochastic equilibrium. It only converges to the most fit allele in the extreme case $p=1$.

Our model is a variation of the well-known statistical physics Ehrenfest chain.
The original Ehrenfest model has no parameter $p$, the transition probabilities are 
$$p(k,k+1)=\frac{N-k}{N}\mbox{ and } p(k,k-1)=\frac{k}{N},$$
see Bhattacharya and Waymire (1990). Our model is a particular case of a two parameter Ehrenfest chain, see Krafft and Schaefer (1993). There, state $0$ flips to state 1 with probability $p_0$ and state $1$ flips it to state 0 with probability $p_1$. The stationary distribution $\nu_N$ for this two parameter model is still a binomial with parameters $N$ and $p_0/(p_0+p_1)$. 

Berger and Cerf (2018) study a model similar to ours but where all sites behave independently of each other. At any given time sites flip or not independently of each other. In particular, there may be multiple flips at any given time. Their model is neutral in the sense that  flips from $0$ to $1$ and from $1$ to $0$ have the same probability. Interestingly, their model also has an Ehrenfest type behavior.

\vskip1cm

{\header The spatial model}

\vskip1cm

We number the sites of the genome from $1$ to $N$.
Let $X_n$ be the configuration of the genome at time $n\geq 0$. That is, for any site $1\leq s\leq N$, $X_n(s)=0$ or $X_n(s)=1$. Note that we can only flip one site at a time in  our model and so $X_n$ and $X_{n+1}$ differ on one site at most. The Markov chain $(Y_n)$ that we studied above is related to the Markov chain $(X_n)$ in the following way. For all $n\geq 0$,
$$Y_n=\sum_{s=1}^N X_n(s).$$
That is, $Y_n$ counts the number of 1's in $X_n$.

Note that $(X_n)$ is a finite irreducible Markov chain with $2^N$ possible states. Therefore, there exists a unique stationary distribution $\pi_N$. For $\eta$ in $\{0,1\}^N$, let
$k(\eta)=\sum_{s=1}^N \eta(s)$ then,
$$\pi_N(\eta)=p^{k(\eta)}(1-p)^{N-k(\eta)}.$$

$\bullet$ The formula for $\pi_N$ shows that as $n\to+\infty$ the probability that a site is in state 1 is $p$ and the probability that a site is in state 0 is $1-p$, independently of all other sites. While the number of sites in state 1 converges to $Np$ their spatial distribution in the genome is completely random.

\medskip

We now show that $\pi_N$ is indeed stationary for $(X_n)$. We will in fact show that it is reversible. Given a configuration $\eta$ a transition is possible to some $\eta'$ if and only if
$\eta(s)=\eta'(s)$ for all $s$ except a single site $s_0$. Assume that $\eta(s_0)=1$ (the case $\eta(s_0)=0$ is treated in a similar way). Then, $\eta'(s_0)=0$. Note that the probability to pick site $s_0$ is $1/N$.
Then, we have the following one step transition probabilities for the Markov chain $(X_n)$,
$$p(\eta,\eta')=\frac{1}{N}(1-p)\mbox{ and }p(\eta',\eta)=\frac{1}{N}p.$$
Using that $k(\eta')=k(\eta)-1$ it is easy to check that 
$$\pi_N(\eta)p(\eta,\eta')=\pi_N(\eta')p(\eta',\eta).$$
Hence, $\pi_N$ is reversible for the chain $(X_n)$.

\vskip1cm

{\header References}

\vskip1cm

R.N. Bhattacharya and E.C. Waymire (1990) {\sl Stochastic processes with applications}, Wiley.

R. Cerf and M. Berger (2018) A basic model of mutations. https://arxiv.org/abs/1806.01212

J.H. Gillepsie (1983) A simple stochastic gene substitution model. Theoretical  Population Biology 23, 202-215.

P. Hegarty and S. Martinsson (2014) On the existence of accessible paths in various models of fitness landscapes. Annals of Applied Probability 24, 1375-1395.

S. Kaufman and S. Levin (1987) Towards a general theory of adaptive walks in rugged landscapes. Journal of theoretical biology 128, 11-45.

H.A. Orr (2003) A minimum on the mean number of steps taken in adaptive walks. Journal of theoretical biology 220, 241-247.

S.C. Port (1994) {\sl Theoretical probability for applications}, Wiley.

R.B. Schinazi (2019) Can evolution paths be explained by chance alone? Journal of Theoretical Biology 465, 65-67.

\end{document}